\input amstex
\documentstyle{amsppt}

\NoRunningHeads
\magnification=\magstep1

\define\e{\varepsilon}

\define\th{\theta}

\define\la{\lambda}
\define\om{\omega}

\define\conv{\text{\rm conv}}
\define\vo #1{\text{\rm vol} \left ( #1 \right )}
\define\etc{, \dots ,}

\define\rn{$\Bbb R^n \,$}
\define\Rn{\Bbb R^n}

\define\nor #1{\left \| #1 \right \|}
\define\enor #1{\Bbb E \, \nor{#1}}
\define\tens #1{#1 \otimes #1}
\define\pr#1#2{\langle {#1} , {#2} \rangle}

\topmatter

\title
Distances between non--symmetric convex bodies and the $MM^*$-estimate.
\endtitle

\author
M. Rudelson
\endauthor

\thanks
This research was  supported in part  by NSF grant DMS-9706835.
\endthanks

\affil
University of Missouri -- Columbia
\endaffil

\address
Department of Mathematics, \hfill \break
University of Missouri, \hfill \break
Mathematical Sciences building, \hfill \break
Columbia, MO 65211
\endaddress

\email
Rudelson\@math.missouri.edu
\endemail

\vskip 1in

\abstract
Let $K, D$ be $n$-dimensional convex bodes. 
Define the distance between $K$ and $D$ as 
$$
d(K,D) = \inf \{ \la \ | \ T K \subset D+x \subset \la \cdot TK \},
$$
where the infimum is taken over all $x \in \Rn$ and all 
invertible linear operators $T$.
Assume that $0$ is an interior point of $K$ and define
$$
M(K) =\int_{S^{n-1}} \nor{\om}_K d \mu (\om), 
$$
where $\mu$ is the uniform measure on the sphere.
We use the difference body estimate to prove that $K$ can be 
embedded into \rn so that 
$$
M(K) \cdot M(K^{\circ}) \le C n^{1/3} \log^a n
$$
for some absolute constants $C$ and $a$.
We apply this result to show that the distance between two $n$-dimensional
 convex bodies does not exceed $n^{4/3}$ up to a logarithmic factor.
\endabstract

\endtopmatter

\vskip .2in

\document

\head 1. Introduction. \endhead

The question of estimating the Banach -- Mazur distance  between
two $n$-dimen\-sional convex symmetric  bodies (i.e. $n$-dimensional
Banach spaces) is one of the central in the Local Theory.
The upper estimate follows from a theorem of John: the distance between 
any convex symmetric body and the ellipsoid does not exceed $\sqrt{n}$.
So, the distance between two such bodies is bounded by $n$.
In 1981 Gluskin \cite{Gl} proved that this estimate is essentially exact.
More precisely, let $m \ge cn$ and let $g_1(\om) \etc g_m (\om)$ be 
independent Gaussian vectors in \rn.
Put 
$$
\Gamma (\om) = \text{abs conv } (g_1(\om) \etc g_m(\om)).
$$
Then with probability close to $1$
$$
d(\Gamma (\om_1),\Gamma (\om_2)) \ge cn.
$$
Here and later $C,c$ etc. mean absolute constants whose value may change from 
line to line.

The situation is entirely different if one considers general convex bodies,
which are not necessary symmetric.
Since for such bodies the origin plays no special role, the definition of the 
distance should be modified to allow shifts.
Namely, the distance between $n$-dimensional convex bodies $K$ and $D$ is
$$
d(K,D) = \inf \{ \la \ | \ T K \subset D+x \subset \la \cdot TK \},
$$
where the infimum is taken over all $x \in \Rn$ and all 
invertible linear operators $T$.

By the theorem of John the distance between a convex body and an ellipsoid 
does not exceed $n$.
This estimate is exact: the distance between the simplex and the ball is
exactly equal to $n$ (actually the simplex is the only convex body having this
 property \cite{P}).
So applying the John theorem twice we show that the distance between two 
convex bodies is bounded by $n^2$.
The lower estimate is very far from this bound.
No examples of non--symmetric bodies the distance between which is greater than
$Cn$ are known.
Moreover, Lassak \cite{L} proved that if 
one of the bodies is symmetric then the distance is bounded by $2n-1$,
while  if one of the bodies is a simplex, it is bounded by $n+2$.
For  $m \ge cn$ one can define a non--symmetric analog of a Gluskin body by
$$
K(\om) = \conv (g_1(\om) \etc g_m(\om)).
$$
Then with probability close to $1$
$$
K(\om) \supset \frac{c}{\sqrt{n}} B_2^n.
$$
This means that the distance between two such bodies is less than $c^2 n$.

It will be shown below that the distance estimate is related to the  
$MM^*$-estimate.
To formulate this more precisely we have to introduce some notation.
Denote by $d_K$  the Banach--Mazur distance between $K$ and 
the Euclidean ball.
Let $\mu$ be the normalized uniform measure on the sphere and let
$$
M(K) = \int_{S^{n-1}} \nor{\om}_K d \mu (\om), \qquad M^*(K)=M(K^{\circ}).
$$
Here $\nor{\cdot}_K$ means the Minkowski functional of $K$ and 
$K^{\circ}=\{ x \in \Rn \ | \ \pr{x}{y} \le 1 \text{ for all } y \in K \}$
is the polar of $K$.
Denote by $g$ a standard Gaussian vector in \rn.
Define the $\ell$-functional by
$$
\ell (K) = \enor{g}_K.
$$
An easy and well known computation shows that
$$
 \frac{1}{\sqrt{n}} \cdot \ell (K) \le  M(K) \le 
\beta_n \frac{1}{\sqrt{n}} \cdot \ell (K),
$$
where $\beta_n \to 1$ when $n \to \infty$.
Denote also by $K_x$ a shift of $K$:
$$
K_x=K-x.
$$

For a convex symmetric body $B$ consider
$$
R(B)=\inf_T M(TB) \cdot M^*(TB),
$$
where the infimum is taken over all invertible {\it affine} 
transformations $T$.
This quantity plays a fundamental role in the Local Theory \cite{M-S},
\cite{P}.
In particular, proofs of the Quotient Subspace Theorem, Inverse Santalo and 
Brunn--Minkowski inequalities are based upon the following estimate 
\cite{P, p.20}
$$
R(B) \le C \log d_B.
$$
This estimate was proved by Pisier.
His proof uses a result of Figiel and Tomczak--Jaegermann \cite{F-T-J}, 
stating that $R(B)$ is controlled by the norm of the Rademacher projection
in $L_2(\Rn, \nor{\cdot}_B)$.

This approach cannot be generalized to the non--symmetric convex bodies,
since the norm of the Rademacher projection can be much bigger.
For example for an $n$-dimensional simplex $\Delta_n$ this norm is 
at least $n/2$, while
$$
c \log n \le R(\Delta_N) \le C \log n.
$$
Using a modified definition of the Rademacher projection, Banasczyk,
Litvak, Pajor and Szarek \cite{B-L-P-S} proved that for any 
$n$-dimensional convex body $K$
$$
R(K) \le C \sqrt{n}.
$$

We shall consider a different approach here.
We shall inscribe a given body into some convex symmetric body for which we
have a good $MM^*$-estimate and compare $M$ and $M^*$ using the comparison 
of the volumes.
Using the volume estimates for the sections of the difference body \cite{R}
we obtain the following results.
\proclaim{Theorem 1}
Any $n$-dimensional convex body $K$ may be embedded in \rn, so that for every 
$\e >0$ there exists a subspace $E \subset \Rn$ of dimension 
at least $(1-\e )n$ such that
$$
\align 
&M^*(K) \le 1 
\intertext{and}
&M(K \cap E) \prec \frac{1}{\e} \cdot \log^2 d_K.
\endalign
$$
\endproclaim

We will write $X \prec Y$ if there exist absolute constants $C, a$ such that
$X \le C \cdot Y \cdot  (\log Y)^a$.
\demo{Remarks}
\item{1.} The Quotient Subspace theorem for convex symmetric 
bodies is based on the $MM^*$-estimate. 
However, the proof of it shows that it is enough to have an  $MM^*$-estimate 
only for a subspace of a small codimension, as in Theorem 1.
\item{2.} An estimate similar to Theorem 1 was also proved by Litvak and 
Tomczak-Jaegermann \cite{L-T-J} by a different method.
Namely, they proved that under the assumptions of Theorem 1 there exists a
subspace $E$ such that 
$$
  M(K\cap E)\, M(K^{\circ}) \le C(\e)  R(K-K) \cdot R(K \cap -K)
  \le  C'(\e)  \log^2{(1+ d_K )},
$$
where $C(\e)$ and $C'(\e)$ depend on $\e$ only.
\enddemo

\proclaim{Theorem 2}
Any $n$-dimensional convex body $K$ may be embedded in \rn, so that
$$
M(K)M^*(K) \prec n^{1/3}.
$$
\endproclaim

The proof of Theorems 1 and 2 consists of several steps.
First, we show in Section 2 that the symmetrization of a convex body
by taking it absolutely convex hull does not affect its volume significantly.
As a byproduct of this observation we obtain the existence of an $M$-ellipsoid
for general convex bodies.
Then in Section 3 we prove that two bodies whose volumes are close possess
a section on which their $\ell$-functionals are close.

In Section 4 we use this fact to prove that the $MM^*$-problem can be
reduced to the estimate of the volumes of sections of the difference body.
By the difference body of a given convex body $K$ we mean
$$
K-K = \{ x-y \ | \ x,y \in K \}.
$$
Then to prove Theorem 1 we apply the following result \cite{R}.
\proclaim{Theorem 3}
Let $K \subset \Rn$ be a convex body and let $F \subset \Rn$ be 
an $m$-dimensional subspace. 
Then
$$
\vo{(K-K) \cap F} \le C^m  \varphi^m (m,n) \cdot  \sup_{x \in \Rn}
\vo{K \cap (F+x)},
$$
where
$$
\varphi (m,n) = \min \left ( \frac{n}{m}, \sqrt{m} \right ).
$$
\endproclaim

Since $\varphi(m,n) \le n^{1/3}$, we have the following immediate
Corollary, which will be used to prove Theorem 2.
\proclaim{Corollary 1}
Let $K \subset \Rn$ be a convex body and let $F \subset \Rn$ be 
an $m$-dimensional subspace. 
Then
$$
\vo{(K-K) \cap F} \le \left (C \cdot  n^{1/3} \right )^m \cdot  
\sup_{x \in \Rn} \vo{K \cap (F+x)}. 
$$
\endproclaim

Also in Section 4 we derive from Theorem 1 the Quotient Subspace Theorem 
for general covex bodies. 
More precisely, applying Theorem 1 and an iteration argument similar 
to \cite{P, Chapter 8}, we obtain the following
\proclaim{Theorem 4}
Let $K$ be a convex body in \rn.
Then for any $\e >0$ there exist linear subspaces $E_2 \subset E_1 \subset \Rn$
such that $\text{\rm dim }E_2 > (1-\e) n$ and for a body 
$D= P_{E_2} (K \cap E_1)$ one has 
$$
d_D \prec \frac{1}{\e^2}.
$$
\endproclaim
\demo{Remark}
It follows from the proof of Theorem 4 that the position of the origin is 
not important here.
In particular, one can assume that the origin coincides with the center of
mass of $K$.
\enddemo

In Section 5 we show that the distance between two convex bodies can be 
estimated in terms of $MM^*$.
More precisely, we prove 
\proclaim{Theorem 5}
Let $K$ and $D$ be $n$-dimensional convex bodies.
Then
$$
d(K,D) \le C n \cdot \sqrt{\log n} \cdot \big ( M(K)M^*(K)+M(D)M^*(D) \big ).
$$ 
\endproclaim
Combinig Theorems 2 and 5, we obtain the following
\proclaim{Corollary 2}
Let $K$ and $D$ be $n$-dimensional convex bodies.
Then
$$
d(K,D) \prec n^{4/3}.
$$
\endproclaim

\demo{Acknowledgment}
I would like to thank E.~Gluskin and A.~Litvak for helpful remarks.
\enddemo

\head 2. Volume estimates. \endhead

Recall that for any convex body $K$ there is a unique point $s \in K$
such that for any point $z$ in the interior of $K$
$$
\vo{K_x^{\circ}} \le \vo{K_z^{\circ}}.
$$
This point is called the Santalo point of $K$.

\proclaim{Lemma 1}
Let $K$ be an $n$-dimensional convex body and suppose that 0 is the Santalo 
point of $K$. Then
$$
\left ( \frac{\vo{\conv \,  (K,-K)}}{\vo{K \cap (-K)}} \right )^{1/n} \le C.
$$
\endproclaim

\demo{Proof}
Using a difference body inequality due to  Rogers and Shephard \cite{R-S},
we obtain that
$$
\left ( \frac{\vo{\conv \,  (K,-K)}}{\vo{B_2^n}} \right )^{1/n} \le
\left ( \frac{\vo{K-K}}{\vo{B_2^n}} \right )^{1/n}\le
C \cdot \left ( \frac{\vo{K}}{\vo{B_2^n}} \right )^{1/n}.
$$
Now applying consecutively  Santalo \cite{S, p.~421} Rogers--Shephard  and 
inverse Santalo \cite{P, p.~100} inequalities, we show that the 
quantity above does not exceed
$$
C \cdot \left ( \frac{\vo{K^{\circ}}}{\vo{B_2^n}} \right )^{-1/n} \le
C \cdot \left ( \frac{\vo{K^{\circ}+(- K^{\circ}}}{\vo{B_2^n}} \right )^{-1/n}
\le C \cdot \left ( \frac{\vo{K \cap (-K)}}{\vo{B_2^n}} \right )^{1/n}.
$$
\qed
\enddemo

The existence of an $M$-ellipsoid for a convex body immediately follows from
Lemma 1. More precisely, for any convex body $K \subset \Rn$ there exists an 
ellipsoid $\Cal E \subset \Rn$ such that
$$
\left ( \frac{\vo{K+ \Cal E}}{\vo{K \cap \Cal E}} \cdot
\frac{\vo{K^{\circ}+ \Cal E^{\circ}}}
{\vo{K^{\circ} \cap \Cal E^{\circ}}}
\right )^{1/n} \le C
$$
This remarkable result was proved by Milman for symmetric bodies.
If $K$ is non-symmetric, assume that $0$ is the Santalo point of $K$ and 
let $\Cal E$ be an $M$-ellipsoid of $D=K \cap (-K)$.
Then Lemma 1 implies that the covering number 
$$
N(K-K,D) \le 2^n \frac{\vo{K-K}}{\vo{D}} \le C^n.
$$
So,
$$
\align
\left ( \frac{\vo{K+ \Cal E}}{\vo{K \cap \Cal E}}\right )^{1/n} &\le
\left ( \frac{\vo{(K-K)+ \Cal E}}{\vo{K \cap (-K) \cap \Cal E}}\right )^{1/n}\\
&\le \left ( \frac{N((K-K),D) \cdot \vo{(K \cap (-K)) + \Cal E}}
{\vo{K \cap (-K) \cap \Cal E}}
\right )^{1/n} \le C.
\endalign
$$
By the inverse Santalo inequality, we have
$$
\align
\left ( \frac{\vo{K^{\circ}+ \Cal E^{\circ}}}
{\vo{K^{\circ} \cap \Cal E^{\circ}}} \right )^{1/n}  
&\le \left ( \frac{\vo{(K \cap (-K))^{\circ} + \Cal E^{\circ}}}
{\vo{(K -K)^{\circ}  \cap \Cal E^{\circ}}}
\right )^{1/n} \\
&\le C
\left ( \frac{\vo{(K-K)+ \Cal E}}{\vo{K \cap (-K) \cap \Cal E}}\right )^{1/n}
\le C.
\endalign
$$
\demo{Remark}
Lemma 1 and the existence of the $M$-ellipsoid were independently proved by
Milman and Pajor \cite{M-P}.
\enddemo

\head 3. $M$-estimates for small sections. \endhead

To prove Theorems 1 and 2 we shall compare $M$ and $M^*$ of a general
convex body with those of certain symmetric bodies.
First notice that if 
$$
B= K-K, \qquad D= K \cap (-K),
$$
then
$$
\align
&\frac12 M^*(B) \le M^*(K) \le M^*(B).
\intertext{By duality we have}
&\frac12 M(D) \le M(K) \le M(D).
\endalign
$$
We shall linearly embed $B$ and thus $K$ into \rn so that
$$
M(B)M^*(B) \le R(B).
$$
Now it is enough to compare $M(B)$ with the minimum of $M(K_x)$ over
all shifts.

By Lemma 1 the volume of $D$ is of the same order as that of $B$.
However the volume estimate alone is not enough to conclude that
$M(D)$ and $M(B)$ are of the same order.
Indeed,  consider the following convex symmetric bodies 
$B,D \subset \Bbb R^{n+1}$:
$$
B=B_2^n+[-e_{n+1},e_{n+1} ], \qquad 
D=B_2^n + [-2^{-n} e_{n+1}, 2^{-n} e_{n+1} ].
$$
Then 
$$
\left ( \frac{\vo{B}}{\vo{D}} \right )^{1/(n+1)} = 2^{n/(n+1)} <2.
$$
However,
$$
M(B) \le 1,
$$
while
$$
M(D) = \frac{1}{\sqrt{n}} \cdot \enor{g}_D \ge 
\frac{1}{\sqrt{n}} \cdot \enor{\gamma e_{n+1}}_D = 
\sqrt{\frac{2}{\pi}} \cdot \frac{2^n}{\sqrt{n}},
$$
where $\gamma$ is a normal random variable.

Nevertheless, it turns out that two convex symmetric bodies of approximately 
the same volume possess a large section on which their $M$-s are close
enough.
More precisely, we prove the following

\proclaim{Lemma 2}
Let $D \subset B$ be $m$-dimensional convex symmetric bodies, such that
$$
\left ( \frac{\vo{B}}{\vo{D}} \right )^{1/m} \le A.
$$
Then  for any $a<1$ there exists a subspace 
$E \subset \Bbb R^m, \ \text{dim} E \ge am$, such that
$$
M(D \cap E) \le \frac{C}{\sqrt{a}} \cdot R(D) \cdot A \cdot M(B) \cdot
 \Big ( CA \cdot M(B) \cdot M^*(B) \Big )^{a/(1-a)}.
$$
\endproclaim
\demo{Remark}
Later we shall choose $a$ so that 
$\Big ( CA \cdot M(B) \cdot M^*(B) \Big )^{a/(1-a)} \le C$.
\enddemo
\demo{Proof}
Let $S$ be an operator such that 
$$
M(SD) \le R(D), \qquad M^*(SD)=1.
$$
Without loss of generality, we may assume that $S$ is self--adjoint.
Let $T=id+M^*(D) \cdot S$.
Then 
$$
\align
&M(TD) \le \frac{M(SD)}{M^*(D)} \le \frac{R(D)}{M^*(D)}, \tag 1\\
&M^*(TD) \le M^*(D) + M^*(D) \cdot M^*(SD) =2 M^*(D).
\endalign
$$
By the inequality of Urysohn \cite{P, p.6}and the inverse Santalo 
inequality, we get
$$
\align
& \left ( \frac{\vo{TD}}{\vo{B_2^m}} \right )^{1/m} \le M^*(TD) \\
& \left ( \frac{\vo{B}}{\vo{B_2^m}} \right )^{-1/m} \le 
c \cdot \left ( \frac{\vo{B^{\circ}}}{\vo{B_2^m}} \right )^{1/m} \le 
c \cdot M(B).
\endalign
$$
Hence,
$$
\align
( \det T)^{1/m} =& \left ( \frac{\vo{TD}}{\vo{D}} \right )^{1/m} =
\left ( \frac{\vo{TD}}{\vo{B_2^m}} \cdot \frac{\vo{B_2^m}}{\vo{B}} \cdot 
\frac{\vo{B}}{\vo{D}} \right )^{1/m} \le \\
& 2M^*(D) \cdot M(B) \cdot A.
\endalign
$$
Since all the eigenvalues of $T$ are greater than 1, there exists an 
invariant subspace $E$ of $T$ of dimension at least $am$, such that 
$$
\nor{T |_E} \le \Big (2M^*(D) \cdot M(B) \cdot A  \Big )^{\frac{1}{1-a}}.
\tag 2
$$
Recall that $g_E$ be the standard Gaussian vector in $E$. We have
$$
\align
M(D \cap E) \le &\frac{1}{\sqrt{\text{dim} E}} \cdot \enor{g_E}_D =
\frac{1}{\sqrt{\text{dim} E}} \cdot \enor{T |_E \, g_E}_{TD}  \\
\le &\nor{T |_E} \cdot \frac{1}{\sqrt{\text{dim} E}} \cdot \enor{g}_{TD} \le
\frac{\sqrt{m}}{\sqrt{\text{dim} E}} \cdot \nor{T |_E} \cdot M(TD).
\endalign
$$
Then  (1) and (2) imply that
$$
\align
M(D \cap E) \le &\frac{1}{\sqrt{a}} \cdot 
\Big (2M^*(D) \cdot M(B) \cdot A  \Big )^{\frac{1}{1-a}} \cdot 
\frac{R(D)}{M^*(D)} \\
\le &\frac{C}{\sqrt{a}} \cdot R(D) \cdot A \cdot M(B) \cdot
 \Big ( CA \cdot M(B) \cdot M^*(D) \Big )^{a/(1-a)}. 
\endalign
$$
The statement of Lemma 2 follows now from $M^*(D) \le M^*(B)$. \qed
\enddemo
Later it will be more convenient to use $\ell(K)$ and $\ell(K^{\circ})$ 
instead of $M(K)$ and $M^*(K)$.
With these parameters the statement of Lemma 2 can be reformulated as
follows:
$$
\ell (D \cap E) \le C \cdot R(D) \cdot A \cdot \ell (B) \cdot
 \Big ( CA \cdot \ell (B) \cdot \ell (B^{\circ}) /m \Big )^{a/(1-a)}.  \tag 3
$$

\proclaim{Lemma 3}
Let $Q$ be a convex body and let $B$  be convex symmetric body
in $\Bbb R^m$. 
Assume that $Q \subset B$ and
$$
\left ( \frac{\vo{B}}{\vo{Q}} \right )^{1/m} \le A.
$$
Then there exists a subspace $F \subset \Bbb R^m$,
$$
\text{dim }F \ge \frac{m}{\log (A \cdot \ell (B) \cdot \ell(B^{\circ})/m)}
$$
and a shift $Q_x=Q-x$ such that
$$
\ell(Q_x \cap F) \le C A \cdot \log d_Q \cdot \ell(B).
$$
\endproclaim
\demo{Proof}
Let $y \in Q$ be a point such that the minimum for the distance from 
$Q-u$ to  an ellipsoid with the center at the origin is attained for $u=y$
and let $z$ be a Santalo point of $Q$.
Put $x=(y+z)/2$ and denote  $D=Q_x \cap (-Q_x)$. Put $x=(y+z)/2$.
Then 
$$
\align
d_{Q_x \cap (-Q_x)} &\le 3 d_Q
\intertext{and by Lemma 1}
\left ( \frac{\vo{B}}{\vo{Q_x \cap (-Q_x)}} \right )^{1/m} &\le 
\left ( \frac{\vo{B}}{\vo{Q}} \right )^{1/m} \cdot 
\left ( \frac{\vo{Q}}{\vo{Q_x \cap (-Q_x)}} \right )^{1/m} \le CA.
\endalign
$$
Put
$$
a=\frac{1}{\log \Big (A \cdot \ell (B) \cdot \ell (B^{\circ}) /m \Big )}
$$
and apply (3).
Since 
$$
 \Big (A \cdot \ell (B) \cdot \ell (B^{\circ}) /m \Big )^{a/(1-a)} \le C,
$$
we have that 
$$
\ell(Q_x \cap F) \le \ell(D \cap F) \le C R(D) \cdot A \cdot \ell(B).
$$
The statement of the Lemma follows now from a theorem of Pisier, since
$$
R(D) \le C \log d_D \le C \log d_Q. \qed
$$
\enddemo

To prove Theorem 2 we need also the following estimate for small dimensions.
\proclaim{Lemma 4}
Let $K$ be an $n$-dimensional convex body.
For any subspace $F \subset \Rn$ there exists a shift $K_u$ such that
$$
\ell (K_u \cap F) \le 2 (\dim F)^2 \cdot \ell (K-K).
$$
\endproclaim
\demo{Proof}
Let $\dim F =m$ and let $e_1 \etc e_m$ be an orthonormal basis of $F$.
Denote $F_j = \text{span } (e_j)$.
We have that
$$
\ell ((K-K) \cap F) = \enor{P_F g}_{K-K} \le
\sum_{j=1}^m  \enor{P_{F_j} g}_{K-K} = 
\sqrt{\frac{2}{\pi}} \sum_{j=1}^m  \nor{e_j}_{K-K}.
$$
Notice that $\nor{e_j}_{K-K} = 1/ \max \nor{x-y}_2$, where the maximum 
is taken over all $x,y \in K$ such that $x-y$ is parallel to $e_j$.
Suppose that the maximum above is attained for the points $x_j, y_j \in K$ 
and put $u_j = 1/2 (x_j+ y_j)$.
Then
$$
\nor{e_j}_{K-u_j} \le \frac{1}{\nor{x_j-u_j}_2} 
=2 \cdot \frac{1}{\nor{x_j-y_j}_2} = 2 \cdot \nor{e_j}_{K-K}.
$$
Let $u=1/m \sum_{j=1}^m u_j$.
Then $K_u \cap F_j \supset 1/m \, (K_{u_j} \cap F_j)$, so 
$$
\align
\ell (K_u \cap F) 
&\le \sum_{j=1}^m \ell (K_u \cap F_j) 
\le m \sum_{j=1}^m \ell (K_{u_j} \cap F_j) 
= m \sum_{j=1}^m \sqrt{\frac{2}{\pi}} \nor{e_j}_{K-u_j} \\
&\le 2m \sum_{j=1}^m \sqrt{\frac{2}{\pi}} \nor{e_j}_{K-K}
\le 2m^2 \cdot \ell ((K-K) \cap F).                          \qed
\endalign
$$
\enddemo

\head 4. Proof of the $MM^*$ estimates.    \endhead
Assume that the body $K$ is embedded into \rn so that
$$
\align
\ell(K-K)           &\le R(K-K) \le C \sqrt{n} \log d_{K-K} 
\le C \sqrt{n} \log d_K              \tag 4 \\
\intertext{and}
\ell((K-K)^{\circ}) &\le \sqrt{n}. \tag 5
\endalign
$$
Since $\ell (K_x)= \Bbb E \sup_{y \in K_x} \pr{g}{y}$ is independent of $x$
for $x \in \text{Int } (K)$,
the last inequality means that for any $x$
$$
\ell(K_x^{\circ}) \le  \sqrt{n}.
$$

To prove theorems 1 and 2 we combine Theorem 3 with the following
\proclaim{Proposition}
Let $K$ be an $n$-dimensional convex body and let $\e> c \log n /n$.
Let $A \le n$ and
assume that for any linear subspace $F \subset \Rn,\ \dim F =m \ge \e n$
there exists a shift $K_x$ such that
$$
\left ( \frac{\vo{(K-K) \cap F}}{\vo{K_x \cap F}} \right )^{1/m} \le A.
$$
Suppose that $K$ is embedded into \rn so that the inequalities (4) and (5)
hold.
Then there exists a subspace $E, \ \dim E \ge (1-\e)n$ and a shift $K_u$
such that
$$
\frac{1}{\sqrt{n}} \cdot \ell (K_u \cap E) 
\prec A \cdot \log^2 d_K \cdot \log^4 1/\e.
$$
\endproclaim
\demo{Proof}
Let $t \in (\e,1)$ and define
$$
\varphi (t) =\min_{x,E} \ell (K_x \cap E), \tag 6
$$
where the minimum is taken over all interior points of $K$ and all linear 
subspaces $E$ of dimension at least $n \cdot (1-t)$. 
We prove the following
\proclaim{Claim}
Let
$$
a=\frac{1}{\log(A \log d_K /\e)}
$$
Then
$$
\varphi^{1/2} ((1-a)t) \le 
\varphi^{1/2} (t) + \left ( C \sqrt{n} A \log^2 d_K \right )^{1/2}.
$$
\endproclaim
To prove the Claim, choose $x$ and $E$, dim$E =n\cdot (1-t)$ so that 
the minimum in (6) occurs for them.
By assumption there exists a shift $y \in K$ such that
$$
\left ( \frac{\vo{(K-K) \cap E^{\perp}}}
{\vo{K_y \cap E^{\perp}}} \right )^{1/nt} \le A.
$$
Put $B= (K-K) \cap E^{\perp}, \ Q = K_y \cap E^{\perp}$.
Then
$$
\ell(B) \le \ell(K-K) \le \sqrt{n} \text{\quad  and \quad } 
\ell(B^{\circ}) \le \ell((K-K)^{\circ}) \le C \sqrt{n} \log n.
$$
Applying Lemma 3
we obtain another shift $K_z$ and a subspace $F \subset  E^{\perp}$  such that
$$
\align
\ell(K_z \cap F) 
&\le C A \log d_Q \cdot \ell(B) \le C A \log d_K \cdot \ell(K-K) \\
&\le C A \log^2 d_K \cdot \sqrt{n}.
\endalign
$$
Here
$$
\align
\dim F 
&\ge \frac{\dim E^{\perp}}
{\log (A \cdot \ell (B) \cdot \ell^*(B)/\dim E^{\perp})}
\ge \frac{nt}{\log (A \cdot C n \log d_K /nt)} \\
&\ge \frac{nt}{\log (A \log d_K /\e)} = ant.
\endalign
$$
Since $A \le n$ and $t> \e > c \log n /n$, we have $\dim F \ge 1$.

For $\mu \in (0,1)$ let $u= \mu x + (1-\mu) z$.
Then $K_u \supset \mu K_x$ and $K_u \supset (1-\mu) K_z$, so
$$
\align
\ell (K_u \cap E) &\le \mu^{-1} \cdot l(K_x \cap E) = \mu^{-1} \cdot \phi (t)
\intertext{and   }
\ell (K_u \cap F) &\le (1-\mu)^{-1} l(K_x \cap F)\le
(1-\mu)^{-1} \cdot CA \log^2 d_K \cdot \sqrt{n}.
\intertext{Since }
\ell(K_u \cap (E \oplus F)) &\le \ell(K_u \cap E) + \ell(K_u \cap F)
\endalign
$$
and $\text{dim }(E \oplus F) \ge n(1-t)+atn=n \cdot (1-(1-a)t)$, we have
that for any $\mu \in (0,1)$
$$
\varphi ((1-a)t) \le 
\mu^{-1} \varphi (t) + (1-\mu)^{-1} \cdot CA \log^2 d_K \cdot \sqrt{n}.
$$

Denote
$$
v= \varphi (t), \qquad w= CA \log^2 d_K \cdot \sqrt{n}. \tag 7
$$
Choose now $\mu \in (0,1)$ so that
$$
\th(\mu) = \mu^{-1} v + (1-\mu)^{-1} w
$$
will be minimal. 
For $\mu_0 = (1-\sqrt{w/v})/ (1-w/v)$ we have
$$
\th(\mu) = (\sqrt{v} + \sqrt{w})^2. \tag 8
$$
Substituting (7) into (8), we get the statement of the Claim. \qed

Now we complete the proof of the Proposition.
Iterating the Claim we obtain that for any $l$
$$
\varphi^{1/2}((1-a)^l) \le 
\varphi^{1/2}(1) + l \cdot (CA \log^2 d_K \cdot \sqrt{n})^{1/2},
$$
provided that $(1-a)^l \ge \e$.
Choose $l$ so that $(1-a)^l \ge \e \ge (1-a)^{l+1}$.
Then for some subspace $\tilde E, \dim  \tilde E \ge n (1-\e)$ and for 
some shift $K_v$ we have 
$$
\ell(K_v \cap \tilde E) \le \varphi((1-a)^{l+1}) \le 
C \cdot l^2 \cdot (CA \log^2 d_K \cdot \sqrt{n}).              
$$
Since 
$$
l \le c \frac{\log 1/\e}{\log (1-a)} \le Ca^{-1} \log 1/\e \le
C \log (A \log d_K) \cdot \log^2 1/\e,
$$ 
we have that
$$
\frac{1}{\sqrt{n}} \cdot \ell(K_v \cap \tilde E) 
\prec A \cdot \log d_K \cdot \log^4  1/\e.   \qed
$$  
\enddemo

\demo{Proof of Theorem 1}
We shall combine Theorem 3 and the Proposition. 
By Theorem 3 for any $m \le n$ and any $m$-dimensional subspace $F$ there 
exists a shift $K_x$ such that
$$
\left ( \frac{\vo{(K-K) \cap F}}{\vo{K_x \cap F}} \right )^{1/m} 
\le C \frac{n}{m}.
$$
So, by the Proposition, there exists a subspace $E, \ \dim E \ge (1-\e)n$
and a shift $K_u$ such that
$$
\frac{1}{\sqrt{n}} \ell (K_u \cap E) 
\prec \e^{-1} \cdot \log^2 d_K \cdot \log^4 \e^{-1} 
\prec \e^{-1} \cdot \log^2 d_K.
$$
We shall show that the point $u$ can be chosen independently of the dimension.

For $1\le l \le L=\log_2 n$ put $\e_l = 2^{-l}$.
Obviously, it is enough to prove Theorem 1 for $\e = \e_l, l= 1 \etc L$. 
By the Proposition we can find a shift $K_{u_l}$ and
a subspace $E_l, \ \dim E_l \ge (1-2^{-l}) n$ such that
$$
\frac{1}{\sqrt{n}} \ell (K_{u_l} \cap E_l) 
\prec 2^l \cdot \log^2 d_K.
$$
Define 
$$
v=\frac{1}{s} \sum \limits_{l=1}^L \frac{u_l}{l^2},
$$
where $s= \sum_{l=1}^L 1/l^2$.
Then $v \in K$.
The same argument as before yields that for any $l$
$$
\frac{1}{\sqrt{n}} \ell (K_v \cap E_l) 
\prec l^2 \cdot 2^l \cdot \log^2 d_K \prec 2^l \cdot \log^2 d_K. \qed
$$
\enddemo

\demo{Proof of Theorem 2}
We shall use Corollary 1 to estimate the sections of the difference body.
For any linear subspace $F \subset \Rn$ there exists a shift $K_x$
such that 
$$
\left ( \frac{\vo{(K-K) \cap F}}{\vo{K_x \cap F}} \right )^{1/m} 
\le A =Cn^{1/3}.
$$
Applying the Proposition with $\e= C \log n /n$, we can find a 
subspace $E, \ \dim E \ge (1-\e)n$ and a shift $K_u$ such that 
$$
\frac{1}{\sqrt{n}} \ell (K_u \cap E) \prec n^{1/3} \log^6 d_K \prec n^{1/3}.
$$
Let $F=E^{\perp}$. 
Then $\dim F = C \log n$ and by Lemma 4 there exists a shift $K_v$ 
such that
$$
\ell (K_v \cap F) \le (C \log n)^2 \cdot \ell (K-K) 
\le  (C \log n)^2 \cdot C \sqrt{n} \cdot \log d_K \prec \sqrt{n}.
$$
Put $w=(u+v)/2$.
Then 
$$
\ell (K_w) \le \ell (K_w \cap E) + \ell (K_w \cap F)
\le 2 \ell (K_u \cap E) + 2 \ell (K_v \cap F) 
\prec n^{5/6} +n^{1/2} \prec n^{5/6}.
$$
Since 
$$
\ell (K_w^{\circ}) = \ell (K^{\circ}) \le \sqrt{n},
$$
this completes the proof of Theorem 2.   \qed
\enddemo

The proof of Theorem 4 follows closely the standard iteration argument
\cite{P, Chapter 8}, so we shall only sketch it.
\demo{Proof of Theorem 4 (Sketch)}
There exists a linear operator $T: \Rn \to \Rn$ and a point
$x \in \Rn$ such that the body $TK+x$ satisfies the estimates of Theorem 1.
Obviously, it is enough to prove Theorem 4 for $TK$ instead of $K$,
so we shall assume that $T= id$.

Let $Q_x$ be the orthogonal projection with 
$\text{Ker }Q_x = \text{span} (x)$.
Notice that it is enough to construct a projection of a section through 0 
of the body $\bar K= Q_x K$.
For this body we have 
$$
M^*(\bar K) \le 2 M^*(K) \le 2
$$
and for any $\e>0$ there exists a linear subspace $E$ of dimension at least
$(1-\e)n$ such that
$$
M(\bar K \cap E) \prec \frac{1}{\e} \log d_{\bar K}.
$$
It follows from \cite{P, Theorem 5.8} that there exists a linear subspace
$V \subset E$ of dimension bigger than $(1-2\e)n$ such that
$$
\bar K \cap V \subset (\bar K - \bar K) \cap V 
\subset \la \cdot (B_2^n \cap V), \tag 9
$$
where
$$
\la = \frac{C}{\e^{1/2}} M^*(\bar K - \bar K) \le 
\frac{C'}{\e^{1/2}} M^*(\bar K) \le \frac{C''}{\e^{1/2}}.
$$
Applying \cite{P, Theorem 5.8} again, we show that there exists a further
subspace $W \subset V$ such that $\text{dim }W > (1-3\e)n$ and
$$
P_{W} (\bar K \cap V) \supset 
P_{W} ((\bar K \cap (-\bar K)) \cap V) \supset 
\mu^{-1} \cdot (B_2^n \cap W),  \tag 10
$$
where 
$$
\mu = \frac{C}{\e^{1/2}} M(\bar K \cap (- \bar K)) \le 
\frac{C}{\e^{1/2}} M(\bar K)  \prec \frac{1}{\e^{3/2}} \log d_{\bar K}.
$$
Finally, it follows from (9) and (10) that
$$
d_{P_{W}} (\bar K \cap V) \prec \frac{1}{\e^2} \log d_{\bar K}.
$$
Theorem 4 follows now from Milman's iteration argument similar to 
\cite{P, Lemmas 8.5 and 8.6}.     \qed
\enddemo

\head 5. Distances between convex bodies. \endhead

We shall use the following Theorem due to Benyaminy and Gordon \cite{B-G}.
\proclaim{Theorem 6}
Let $K$ and $D$ be $n$-dimensional convex bodies.
then
$$ 
\align
d(K,D) \le \frac{C}{n} &\cdot \Big ( 
\nor{id: B_2^n \to K^{\circ}} \cdot \ell (D) 
+ \nor{id: B_2^n \to D} \cdot \ell (K^{\circ}) \Big ) \\
&\times \Big ( \nor{id: B_2^n \to D^{\circ}} \cdot \ell (K) 
+ \nor{id: B_2^n \to K} \cdot \ell (D^{\circ}) \Big ).
\endalign
$$
\endproclaim
Notice that although Theorem 6 was proved under the assumption that 
the bodies $K$ and $D$ are symmetric, the same proof works for general 
convex bodies.

By a standard contraction argument
$$
\nor{id: B_2^n \to K} \le \ell (K),
$$
so by Theorem 2 we have
$$
d(K,D) \le  \frac{C}{n} \cdot \ell(K^{\circ}) \ell(D) \cdot
\ell(D^{\circ}) \ell(K).
$$
However combining this approach with a result of Banasczyk, Litvak, Pajor
and Szarek \cite{B-L-P-S, Proposition 3.1 and Remark 3.2} (see also \cite{B}),
 we obtain a better estimate.
More precisely, we need the following
\proclaim{Theorem 7}
Let $K$ be an $n$-dimensional convex body and let $B_2^n$ be the 
ellipsoid of maximal volume inscribed in $K$.
Then
$$
\ell (K^{\circ}) \le C n \sqrt{\log n}.
$$
\endproclaim

\demo{Proof of Theorem 5}
Let $W=\max \Big ( M(K)M^*(K), M(D)M^*(D) \Big )$.
It is enough to prove that
$$
d(K,D) \le Cn \cdot \sqrt{\log n} \cdot W.
$$
Assume that the body $K$ is embedded into \rn so that
$$
\ell (K) \le \sqrt{n} \quad \text{end} \quad 
\ell (K^{\circ}) \le \sqrt{n} \cdot W.
$$
Let $S$ be a linear operator which maps the ellipsoid of maximal volume 
inscribed in $K$ onto $B_2^n$.
Put
$$
T=id + \frac{W}{\sqrt{n \log n}} S.
$$
Then
$$
\nor{id: B_2^n \to TK} \le 
\nor{id: B_2^n \to \frac{W}{\sqrt{n \log n}} SK} \le \frac{\sqrt{n \log n}}{W}.
\tag 11
$$
Also we have
$$
\ell (TK) \le \ell (K) \le  \sqrt{n} \tag 12 
$$
and by Theorem 7
$$
\aligned
\ell ((TK)^{\circ}) &\le \ell (K^{\circ}) + 
\frac{W}{\sqrt{n \log n}} \cdot \ell ((SK)^{\circ}) \\
&\le \sqrt{n} \cdot W + \frac{W}{\sqrt{n \log n}} \cdot Cn \sqrt{\log n}
\le C \sqrt{n} \cdot W.    
\endaligned \tag 13
$$
By the contraction argument
$$
\nor{id: B_2^n \to (TK)^{\circ}} \le \ell ((TK)^{\circ}) \le 
C \sqrt{n} \cdot W.        \tag 14
$$
Similarly, there exists an embedding $U$ of $D$ into \rn such that
$$
\align
&\nor{id: B_2^n \to (UD)^{\circ}} \le \frac{\sqrt{n \log n}}{W}, \qquad
\ell ((UD)^{\circ}) \le \sqrt{n} 
 \intertext{and} 
&\nor{id: B_2^n \to UD} \le \ell (UD) \le C \sqrt{n} \cdot W.
\endalign
$$
Substituting these estimates and (11) -- (14) into Theorem 6  we obtain the 
statement of Theorem 5.       \qed

\enddemo

\vskip .5in
\subheading{References}
\medbreak

\item{[B]}
F. Barthe, An extremal property of the mean width of the simplex, Math. Ann.
{\bf 310} (1998), no. 4, 685--693.

\item{[BG]}
Y.~Benyamini, Y.~Gordon,  {\it Random factorization of operators between Banach
spaces,} J. Analyse Math. {\bf 39} (1981), 45-74.

\item{[B-L-P-S]}
W.~Banaszczyk, A.~Litvak, A.~Pajor and S.~Szarek, 
{\it The flatness theorem for
non--symmetric convex bodies via the local theory of Banach spaces}, Preprint.

\item {[Gl]} 
E.~Gluskin,
{\it The diameter of the Minkowski compactum is approximately equal
to $n$}, Funct. Anal., Appl.  {\bf 15} (1981), 72-73. (in Russian).

\item{[L]}
M.~Lassak, {\it Approximation of convex bodies by centrally symmetric bodies},
Geom. Dedicata {\bf 72} (1998), no. 1, 63--68. 

\item{[L-T-J]}
A.~Litvak, N.~Tomczak-Jaegermann, {\it Random aspects of
high-dimensional convex bodies}, preprint.

\item{[M-P]} 
V.~D.~Milman and A.~Pajor, {\it Entropy and asymptotic geometry
of non-symmetric convex bodies}, preprint.

\item{[M-S]}
V.~D.~Milman and G.~Schechtman, {\it Asymptotyc Theory of Finite
Dimensional Normed Spaces}, Lecture Notes in Math, {\bf 1200}, 1986.

\item{[P]} 
G.~Pisier, {\it The volume of convex bodies and Banach space
geometry},  Cambridge Tracts in Mathematics {\bf 94}, 1989.

\item{[Pa]} 
O.~Palmon,
{\it The only convex body with extremal distance from the ball is the simplex},
Israel J. of Math., {\bf 79} (1993).

\item{[R-S]}
C.~A.~Rogers and  G.~C.~Shephard, {\it The difference body of a convex body},
Arch. Math. {\bf 8} (1957), 220--233.

\item{[R]} 
M.~Rudelson, {\it Sections of the difference body},
 to appear in Discrete and Computational Geometry.

\item{[S]} 
R. Schneider, {\it Convex bodies: the Brunn-Minkowski theory},
 Encyclopedia of Mathematics and its Applications, {\bf 44},
 Cambridge University Press, Cambridge, 1993.

\enddocument